\documentclass{article}

\usepackage{arxiv}

\usepackage[utf8]{inputenc} 
\usepackage[T1]{fontenc}    
\usepackage{hyperref}       
\usepackage{url}            
\usepackage{booktabs}       
\usepackage{amsfonts}       
\usepackage{nicefrac}       
\usepackage{microtype}      
\usepackage{lipsum}
\usepackage{algorithm}
\usepackage{algpseudocode}
\usepackage{graphicx}
\usepackage{amsmath,booktabs,latexsym}
\usepackage{mathptmx}
\usepackage{epsfig}
\usepackage{mathptmx}
\graphicspath{ {./images/} }
\newtheorem{definition}{Definition}

\title{Paradoxes Are Not Contradictions: Re‑examining the Third Mathematical Crisis}

\author{
 Yuan H. Y. \\
  \texttt{yuanhaoyan@126.com} \\
}

\begin{document}
\maketitle
\begin{abstract}
Russell's paradox was proposed in the early 20th century to address loopholes in set theory, which directly triggered the third mathematical crisis. This paper proposes and elaborates a perspective distinct from previous studies: paradoxes do not give rise to contradictions; instead, they constitute a Möbius strip-style self-consistent logical structure via self-reference and negation under the logical rules within a system. Gödel's incompleteness theorems indicate that such structures universally exist in formal logical systems. Turing proved the undecidability of the halting problem by first assuming the existence of a halting program and subsequently refuting this assumption through paradox construction, and this paper demonstrates flaws inherent to such proof strategy. Cases of paradoxes within three-valued logical systems are further discussed in this work, where the undecidability of paradoxes is rigorously proven. Finally, inspirations drawn from paradoxes for the real world are explored: two opposing factors can be integrated through the joint mechanism of self-reference and negation. A representative example is the wave-particle duality of light, whose essence may be interpreted as a paradox of waves and particles.
\end{abstract}


\section{Introduction}
At the beginning of the 20th century, British mathematician Bertrand Russell put forward Russell's paradox, which brought about a crisis in set theory founded by Cantor. In the early 1920s, Hilbert presented an ambitious vision of constructing an axiomatic system where all mathematical propositions could be verified or falsified via finite deductive steps while satisfying consistency and completeness \cite{hilbert}. Consistency, also referred to as non-contradiction, requires that no contradictory propositions coexist within an axiomatic system; that is, $p$ and \(\neg p\) cannot hold simultaneously. Completeness demands that every proposition formulated within the system is either provable or refutable. Nevertheless, Kurt Gödel established the incompleteness theorems in 1931. By constructing Gödel numbers and a specific proposition, he proved that first-order predicate logic incorporating elementary number theory contains true yet unprovable propositions \cite{godel}. This result implies that completeness is unattainable. From the perspective of proof procedures, such undecidable propositions share an analogous structural feature with Russell's paradox, namely the combination of self-reference and negation. Based on this characteristic, these propositions can be systematically investigated, rather than being regarded as an unknowable boundary of truth. Two primary approaches have been developed within academia to resolve Russell’s Paradox. The first introduces restrictive axioms to naive set theory so as to avoid the emergence of Russell’s Paradox, examples including the ZF system and hierarchical stratification frameworks \cite{limit}. The second approach attributes the root of the problem to the unrestricted comprehension principle and thus abandons this principle entirely to construct alternative formal systems \cite{new,hubei}. Focusing primarily on Russell's paradox, this paper proposes two viewpoints that diverge from prevailing mainstream arguments: 1. Paradoxes do not correspond to contradictions but constitute a special, valid logical structure. 2. The emergence of paradoxes is a natural phenomenon; accordingly, the premises giving rise to paradoxes cannot be rejected merely to eliminate paradoxes themselves.

\section{Main results}
\subsection{Paradoxes are not contradictions}
Unless otherwise specified, the terminologies used throughout this paper are defined as follows:
\begin{itemize}
\item{Contradiction}: Logical invalidity arising from the violation of the consistency axiom.
\item{Paradox}: Let $A$ and $B$ denote two distinct logical judgments of an identical proposition. $B$ can be deductively derived under the premise $A$, and conversely, $A$ can be deductively derived under the premise $B$.
\item{Russell's paradox}: A set either belongs to itself (\(S \in S\)) or does not belong to itself (\(S \notin S\)). In accordance with the comprehension principle of Cantorian set theory, let \(S_1=\{x:x\in x\}\) be the set consisting of all sets that contain themselves as elements, and let \(S_2=\{x:x\notin x\}\) be the set consisting of all sets that do not contain themselves as elements. The core question is whether \(S_2\) is an element of itself. If \(S_2\in S_2\), then \(S_2\) should not be a member of \(S_2\); if \(S_2\notin S_2\), then \(S_2\) ought to be a member of \(S_2\).
\item{Barber paradox}: A barber claims that he shaves every person who does not shave themselves. A logical dilemma emerges regarding whether the barber ought to shave himself. If the barber shaves himself, he violates his rule and should not shave himself; if the barber does not shave himself, his rule mandates that he must shave himself.
\end{itemize}
Concerning paradoxes, the mainstream research framework adopts the following reasoning paradigm: an initial hypothesis is postulated, a conclusion contradictory to this hypothesis is deduced, and the original hypothesis is consequently discarded. If all plausible hypotheses are eliminated through this procedure, the underlying proposition is deemed undecidable. To rigorously clarify the definition of contradiction, three distinct logical scenarios are analyzed as follows:

(1) If premises $G$ simultaneously entail propositions $A$ and $B$, where $A$ and $B$ form a contradictory pair, the consistency principle is violated. The premise $G$ must therefore be rejected, and only a contradictory relation exists between $A$ and $B$ in this case.

(2) Under the condition that the premise $G$ holds unconditionally, suppose proposition $A$ is true, from which proposition $B$ is deduced. If $B$ contradicts the universal premise $G$, the consistency principle dictates that the hypothetical proposition $A$ must be refuted. Only a causal logical relation connects $A$ and $B$ here, and this structure widely appears in proof by contradiction.

(3) Given the universal premise $G$, supposing $A$ to be true yields $B$, and supposing $B$ to be true yields $A$. In this scenario, both contradictory and bidirectional causal relations hold between $A$ and $B$, which precisely characterizes the logical structure of paradoxes.

The first two scenarios have attained universal consensus within academia. The core controversy lies in the third case: deducing $B$ from $A$ inherently negates $A$, and the same reciprocal negation holds when deducing $A$ from $B$. Notably, $A$ and $B$ never hold simultaneously, which implies no violation of the consistency principle and the absence of genuine logical contradiction. Accordingly, paradoxes are logically valid constructions.
\begin{figure}[htbp]
	\centering
	\includegraphics[width=0.8\columnwidth]{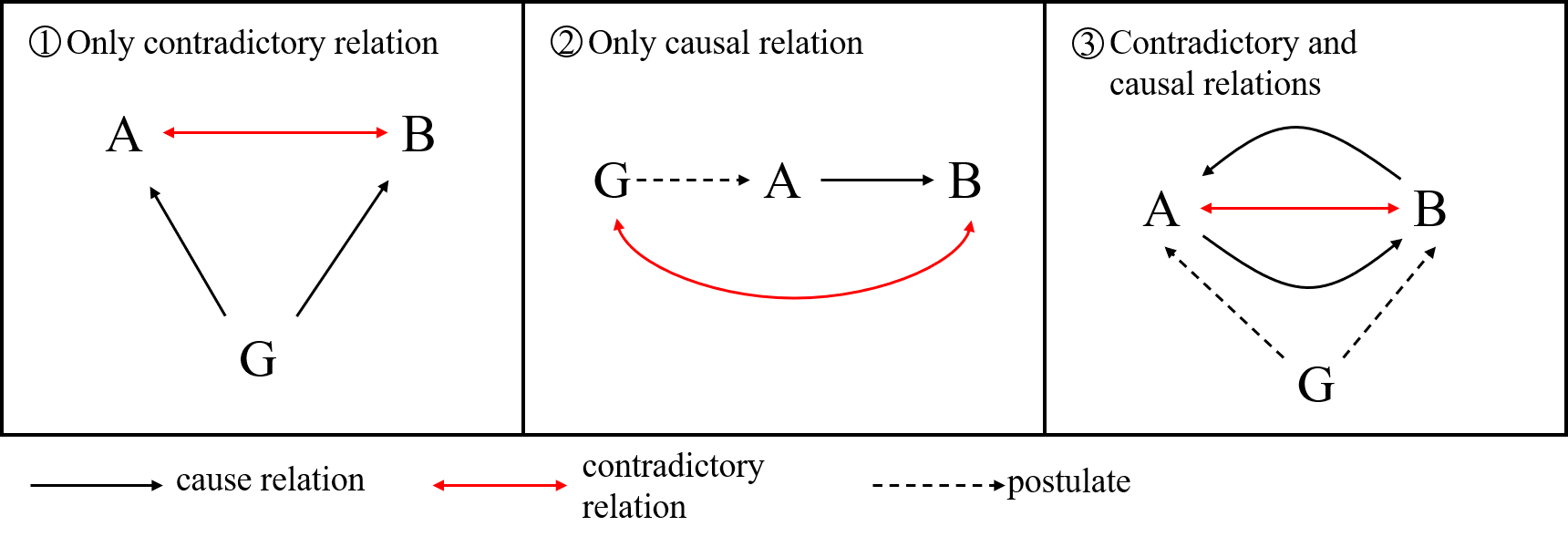}
	\caption{Three cases of contradictions.}
	\label{relation}
\end{figure}

Subsequent rigorous derivations verify that the two propositions $A$ and $B$ involved in Russell's paradox do not violate the consistency principle. According to the definition of \(S_2\), for any arbitrary set x, \(x \notin x\) holds if \(x \in S_2\); if \(x \notin S_2\), then \(x \in S_1\), which equivalently implies \(x \in x\). Substituting \(S_2\) for the variable x, the logical statements are transformed into: "If \(S_2 \in S_2\), then \(S_2 \notin S_2\); if \(S_2 \notin S_2\), then \(S_2 \in S_2\)". This recursive reasoning is essentially equivalent to the logical negation operation. Therefore, the intrinsic nature of a paradox corresponds to a self-negating logical structure constructed by combining negation constraints within a formal system with self-reference. Regardless of whether the self-referential proposition is assumed true or false, a single deductive step yields the opposite truth valuation, meaning the entire deductive procedure acts fundamentally as a negation operation. Conventional logical deductions strictly abide by consistency: valid premises combined with sound inference rules guarantee true conclusions. Since the two propositions before and after a negation operation cannot hold true concurrently, no genuine logical contradiction arises within Russell's paradox.
\subsection{Möbius strip structure}
Given that paradoxes constitute logically valid structures, the following formal definition is presented:
\begin{definition}
Let $A$ and $B$ denote two distinct logical truth valuations of the same proposition. If \(A \rightarrow B\) and \(B \rightarrow A\), the solution to the corresponding problem is defined as a paradox of $A$ and $B$.
\end{definition}
If the truth and falsity of a proposition are regarded as the two sides of a paper strip and logical deduction is interpreted as moving along the strip, a paradox can be analogized to a Möbius strip. Starting from any point on the central line of the strip and moving in a fixed direction, one will inevitably arrive at the reverse side of the starting point before returning to the original position. Suppose a paper strip has a permanently black front side and a permanently white reverse side; let point $A$ (black) be marked at the center of the strip, with its corresponding backside point $B$ (white). When the tail end of the strip is rotated by 180 degrees (negation) and joined to the head end (self-reference), a Möbius strip is formed, and a natural question arises regarding its overall color. If the strip is assumed black, traveling one full length from point $A$ will lead to point $B$, which is white. Symmetrically, moving one full length from point $B$ will return to the black point $A$. Consequently, the Möbius strip can neither be classified as purely black nor purely white, and its color state may be designated as a paradox of black and white. The distinctive property of a Möbius strip lies in its unification of two originally opposing surfaces into a single continuous surface. By analogy, two mutually exclusive logical propositions are unified via the joint mechanism of self-reference and negation. After such unification, standard binary truth assignments of "true" or "false" become inapplicable, and the only valid characterization is the statement of a paradox of true and false, or a paradox of propositions $A$ and $B$.

The analogy of the Möbius strip naturally motivates the conjecture of its symmetric logical counterpart: self-reference without negation, which corresponds to an ordinary ring formed by directly joining the two ends of an untreated paper strip and admits two separate logical truth valuations. Take the membership relation of \(S_1\) in Russell's paradox as an illustrative example. If \(S_1\) is assumed to contain itself, the definition of \(S_1\) directly implies \(S_1\in S_1\); if \(S_1\) is assumed not to contain itself, it follows that \(S_1\in S_2\), which further yields \(S_1\notin S_1\). Regardless of the initial membership assumption for \(S_1\), the derived conclusion consistently corroborates the original hypothesis. Although this reasoning constitutes circular argumentation, it remains logically consistent. In summary, self-referential propositions fall into two symmetric special categories: circular affirmative arguments, where every truth valuation is self-validating, and paradoxes induced by negation, where neither truth valuation can hold.
\section{Paradoxes and the halting problem}
The halting problem constitutes a fundamental problem in computer science, which asks whether a universal procedure can be constructed such that, given any executable program and its corresponding input parameters, the procedure is capable of determining within finite computational steps whether the target program terminates. Turing proved that such a universal procedure cannot exist. The core proof proceeds by supposing that a procedure \(T(P,s)\) with the above property exists, followed by the construction of a self-referential function

\begin{algorithm}

	\begin{algorithmic}[1]
		\State \textbf{B}$(x):\{$
		\State \textbf{If}~ $T(x,x)=1$
		\State ~~~~\textbf{While}$(1);$
		\State \textbf{If}~ $T(x,x)=0$
		\State ~~~~\textbf{Return}~$1;\}$

	\end{algorithmic}
\end{algorithm}

What about whether the program \(B(B)\) halts? If \(T(B,B)\) outputs 1, then \(B(B)\) does not terminate; if \(T(B,B)\) outputs 0, then \(B(B)\) terminates. A seeming contradiction is therefore derived, which is used to conclude that such a universal decision procedure cannot exist.

The two statements that a program halts ($A$) and that a program does not halt ($B$) form a pair of opposing logical truth valuations. The proof establishes \(A \rightarrow B\) and \(B \rightarrow A\), which satisfies the definition of a paradox provided earlier. Nevertheless, Chapter 2 rigorously demonstrates that a paradox represents a logically valid structure rather than a genuine contradiction. Accordingly, the premise assuming the existence of a universal halting decider cannot be rejected merely to eliminate the emergent paradox (see Case 3 in Figure 1). If paradoxes are treated as legitimate logical constructions, the following conclusion can be drawn: any program designed to predict termination behavior may fail when tasked with analyzing another program that contains or invokes the predictor itself.

In fact, paradoxical structures can be constructed not only for termination predictors but for any program that predicts the output of target programs. Following the standard proof logic of the halting problem, any decision procedure capable of forecasting program outputs should be proven non-existent. As an example, consider the following question: can a program be devised that, given any arbitrary program and its input arguments, determines within finite computational steps whether the program crashes with an error message consisting of an extremely long specific key sequence? Intuitively, almost no practical program would generate such an error, which seemingly simplifies the development of the corresponding decision procedure. Analogously, suppose such a procedure \(T'(P,s)\) exists, and a self-referential function is constructed accordingly

\begin{algorithm}
	
	\begin{algorithmic}[1]
		\State \textbf{B'}$(x):\{$
		\State \textbf{If}~ $T'(x,x)=1$
		\State ~~~~\textbf{Return}$(0);$
		\State \textbf{If}~ $T'(x,x)=0$
		\State ~~~~\textbf{Error}~(an extremely long specific key sequence)$;\}$

	\end{algorithmic}
\end{algorithm}

This arise the similar paradox. In other words, for any procedure tasked with predicting the outputs of target programs, regardless of whether its decision task is trivial or intractable, a paradoxical counterexample program can always be constructed to render the predictor ineffective.

Does this fact truly demonstrate that no universal procedure can correctly decide the termination property of all programs? The interpretation presented in this work is that the conventional comprehension of the quantifier "all" is imprecise. Take Russell’s paradox as an illustration: after constructing sets \(S_1\) and \(S_2\), one naturally assumes that all sets can be classified into these two collections. However, mathematical reasoning contradicts this intuition, since \(S_2\) belongs to neither \(S_1\) nor \(S_2\), nor to any other set. This raises the question of whether the comprehension principle is flawed for failing to encompass \(S_2\). The answer is negative—the comprehension principle itself remains sound. Defining \(S_1\) and \(S_2\) inherently generates a paradox: \(S_2\) belongs to the paradox of \(S_1\) and \(S_2\). Every set excluding this paradoxical object can be categorized into either \(S_1\) or \(S_2\), and this restricted scope is the proper interpretation of the term "all".
An analogous reasoning applies to formal logical systems. If only two truth valuations, true and false, are stipulated for well-formed propositions, there inevitably exist statements that resist binary truth assignment—for instance, the liar paradox: "This sentence is false". Nonetheless, the rule restricting propositional truth values to true and false retains its validity. For the halting problem, suppose a procedure exists that correctly decides termination for arbitrary programs; a paradoxical counterexample program can always be constructed on which the decider yields no valid judgment. Meaningful discussion of a universal halting decider is therefore only feasible when such paradoxical programs are excluded from consideration.

In summary, paradoxical programs are logically well-formed, yet the existence of halting deciders cannot be dismissed on this account. In other words, one ought not to reject the premises that give rise to paradoxes merely for the sake of eliminating such paradoxical constructions.

\section{Paradoxes and undecidability}
All previous discussions are grounded in two-valued logical systems, yet paradoxes can also arise within many-valued logics. We illustrate this with a three-valued system via the following formal definition

\begin{definition}
	Let $A$, $B$, and $C$ be three propositions pertaining to the same subject term, pairwise mutually exclusive and collectively exhaustive, jointly forming the complete set of truth valuations for three-valued logic. If
\begin{equation}
\left\{
\begin{aligned}
	A\rightarrow B\\
	B\rightarrow C\\
	C\rightarrow A
\end{aligned}\right.
\end{equation}
The resolution to the corresponding problem is defined as a paradox of $A$, $B$, and $C$.
\end{definition}
Drawing on the earlier Möbius strip analogy, a three-valued paradox can be likened to a triangular prism whose tail end is rotated by $120$ degrees and fused to its head end. Its original three distinct lateral surfaces merge into a single continuous surface, analogous to three pairwise mutually exclusive propositions forming a logical structure via cyclic implication and self-reference.

There exists an alternative scenario where propositions $A$, $B$, $C$ satisfy the system
\begin{equation}
	\left\{
	\begin{aligned}
		A\rightarrow B\\
		B\rightarrow A\\
		C\rightarrow C
	\end{aligned}\right.
\end{equation}
Here A and B form a paradox under self-referential reasoning, while proposition C is individually self-consistent. Take the barber paradox as an example: apart from the two alternatives “the barber shaves himself” and “the barber is shaved by others”, we may introduce a third option “the barber never shaves”. If the barber permanently abstains from shaving, he violates none of his own rules.

It is critical to emphasize that paradoxes do not possess the same categorical status as individual truth valuations A, B. Though we previously described Russell’s paradox as “\(S_2\) belongs to the paradox of \(S_1\) and \(S_2\)”, this paradox itself cannot be defined as a set within the system.

\textbf{Proof}:~
Suppose we attempt to define a set \(S_3\) corresponding to "the paradox of \(S_1\) and \(S_2\)", such that \(S_1\), \(S_2\) and \(S_3\) are pairwise mutually exclusive and collectively exhaustive. Assume \(S_3 \in S_3\). By the definition of \(S_1\), it follows that \(S_3 \in S_1\). If \(S_3 \in S_1\), then \(S_3 \notin S_3\), which by the definition of \(S_2\) implies \(S_3 \in S_2\). From the definition of \(S_3\), we further have \(S_2 \in S_3\), which returns us to \(S_3 \in S_3\). This chain of reasoning yields the system
\begin{equation}
	\left\{
	\begin{aligned}
		S_3 \in S_3\rightarrow S_3 \in S_1\\
		S_3 \in S_1\rightarrow S_3 \in S_2\\
		S_3 \in S_2\rightarrow S_3 \in S_3
	\end{aligned}\right.
\end{equation}
It therefore follows that \(S_3\) belongs to a paradox of \(S_1\), \(S_2\) and \(S_3\). Since \(S_2 \in S_3\), we may further deduce that \(S_2\) belongs to a paradox of \(S_1\), \(S_2\) and \(S_3\), which contradicts the initial premise that \(S_2\) belongs to a paradox of \(S_1\) and \(S_2\). This proves that "the paradox of \(S_1\) and \(S_2\)" cannot itself be treated as a set of the same type as \(S_1\) and \(S_2\).

By parallel reasoning, "the paradox of truth and falsity" is not itself a valid truth valuation. Propositions admit only the two standard truth assignments: true or false. If the statement "'This sentence is false' is a paradox" were assigned the value true, an identical contradictory cycle would be derived. Paradoxes are therefore undecidable. In conclusion, the paradox of $A$ and $B$ does not occupy a categorical level parallel to $A$ and $B$.
\section{Conclusions and insights}
This paper addresses self-referential paradoxes exemplified by Russell’s Paradox and advances two novel viewpoints:
1. A paradox is a Möbius-strip-like logical structure generated by self-reference and negation under the rules of a formal system. It does not introduce genuine inconsistency within the system.
2. One should not reject the foundational premises that give rise to paradoxes merely to eliminate paradoxical constructions.
From this, a flaw exists in Turing’s classical proof of the Halting Problem. Turing’s argument only demonstrates that paradoxical programs inevitably exist within the framework—and such paradoxes ought to be expected rather than used to rule out universal halting deciders. Meaningful analysis of universal termination decision procedures is only viable when these paradoxical counterexamples are excluded from consideration. More generally, discussions of a system’s completeness are only coherent after all paradoxical objects are factored out. Fortunately, human mathematical practice already implicitly adopts this approach: mathematicians routinely employ set theory, mathematical logic and other formal tools while setting paradoxes aside, and no fatal contradictions arise in practical application.

Furthermore, paradoxes yield a tangible philosophical insight: two mutually exclusive opposing properties can be unified via the joint mechanism of self-reference and negation. Wave-particle duality of light serves as a paradigmatic illustration. This work conjectures that the intrinsic nature of light constitutes a paradox of wave and particle states. The particle character of light is directly observable, whereas its wave character cannot be directly measured, which further frames this phenomenon as a paradox of observable versus unobservable properties.

\bibliographystyle{unsrt}  


\end{document}